\documentclass[12pt]{amsart}
\usepackage{amssymb}

\title[On formal maps between generic submanifolds in complex space]{On formal maps between generic submanifolds\\ in complex space}

\author{Jean-Charles Suny\'e}

\address{Universit\'e de Rouen, CNRS, Laboratoire de Math\'ematiques Rapha\"el Salem, Avenue de l'Universit\'e, B.P. 12, 76801 Saint Etienne du Rouvray, France}

\email{jean-charles.sunye@etu.univ-rouen.fr}

\thanks{The author is partially supported by the Amadeus program of the 'Partenariat Hubert Curien'.}


\newtheorem{Thm}{Theorem}[section]
\newtheorem{Cor}[Thm]{Corollary}
\newtheorem{Pro}[Thm]{Proposition}
\newtheorem{Lem}[Thm]{Lemma}

\theoremstyle{definition}
\newtheorem{Def}[Thm]{Definition}

\theoremstyle{remark}
\newtheorem{Rem}[Thm]{Remark}
\newtheorem{Exa}[Thm]{Example}

\numberwithin{equation}{section}

\def\bl{\begin{Lem}}
\def\el{\end{Lem}}
\def\bp{\begin{Pro}}
\def\ep{\end{Pro}}
\def\bt{\begin{Thm}}
\def\et{\end{Thm}}
\def\bc{\begin{Cor}}
\def\ec{\end{Cor}}
\def\bd{\begin{Def}}
\def\ed{\end{Def}}
\def\be{\begin{Exa}}
\def\ee{\end{Exa}}
\def\bpf{\begin{proof}}
\def\epf{\end{proof}}
\def\ben{\begin{enumerate}}
\def\een{\end{enumerate}}
\def\beq{\begin{equation}}
\def\eeq{\end{equation}}
\def\beqar{\begin{eqnarray}}
\def\eeqar{\end{eqnarray}}
\def\brm{\begin{Rem}}
\def\erm{\end{Rem}}

\def\hn {holomorphically nondegenerate} 
\def\hyp{hypersurface}

\def\det{{\text{\rm det}}}

\def\adj{\text{\rm adj }}

\def\Fb{\bar F}
\def\Gb{\bar G}
\def\Hb{\bar H}
\def\zb{\bar z}
\def\wb{\bar w}
\def\Qb{\bar Q}
\def\Hkb{\bar {H}^k}
\def\t{\tau}
\def\z{\zeta}
\def\phi{\varphi}

\def\L{\Lambda}
\def\b{\beta}
\def\a{\alpha}

\def\l{\lambda}
\def\c{\chi}
\def\g{\gamma}
\def\rp{{\rho^\p}}
\def\C{\mathbb C}

\def\N{\mathbb N}

\def\J{\text{\rm Jac}}
\def\d{\partial}

\def\Mp{M^\p}
\def\Qp{{Q^\p}}
\def\Qpb{\bar {Q^\p}}
\def\zp{{z^\p}}
\def\cp{\c^\p}
\def\tp{\t^\p}
\def\wp{{w^\p}}
\def\rpt{{\tilde \rp}}
\def\rpb{{\bar \rp}}

\def\calR{\mathcal R}
\def\Qpa{{\Qp_\a}}
\def\Qpar{{\Qp_\a^r}}
\def\Ub{{\bar U}}
\def\Vb{{\bar V}}
\def\rd{\circ}
\newcommand{\p}{\prime}

\newcommand{\nequiv}{{\equiv \!\!\!\!\!\!  / \,\,}}

\newcommand{\mc}[1]{\mathcal #1}

\newcommand{\deri}[2]{\frac{\d #2}{\d #1}}
\newcommand{\mderi}[3]{\frac{\d^{|#2|} #3}{\d {#1}^{#2}}}
\newcommand{\dblderi}[5]{\frac{\d^{|#2|+|#4|} #5}{\d{#1}^{#2} \d {#3}^{#4}}}
\newcommand{\dblderihyp}[5]{\frac{\d^{|#2|+#4} #5}{\d{#1}^{#2} \d {#3}^{#4}}}

\newcommand{\trpderihyp}[7]{\frac{\d^{|#2|+|#4|+#6} #7}{\d{#1}^{#2} \d{#3}^{#4} \d{#5}^{#6}}}
\newcommand{\dop}[1]{\frac{\d}{\d #1}}
\newcommand{\mdop}[2]{\frac{\d^{|#2|}}{\d {#1}^{#2}}}
\newcommand{\mdopz}[2]{\left.\frac{\d^{|#2|}}{\d {#1}^{#2}}\right|_{#1=0}}
\newcommand{\dbldop}[4]{\frac{\d^{|#2|+|#4|} }{\d {#1}^{#2} \d {#3}^{#4}}}
\newcommand{\dbldophyp}[4]{\frac{\d^{|#2|+#4} }{\d {#1}^{#2} \d {#3}^{#4}}}
\newcommand{\dbldopz}[4]{\left.\dbldop{#1}{#2}{#3}{#4}\right|_{\begin{subarray}{c} 
#1=0\\ 
#3=0
\end{subarray}}}

\begin{document}

\begin{abstract}
Let $H:(M,p)\rightarrow (\Mp,p^\p)$ be a formal mapping between two germs of real-analytic generic submanifolds in $\C^N$ with nonvanishing Jacobian. Assuming $M$ to be minimal at $p$ and $\Mp$ \hn \ at $p^\p$, we prove the convergence of the mapping $H$. As a consequence, we obtain a new convergence result for arbitrary formal maps between real-analytic hypersurfaces when the target does not contain any holomorphic curve. In the case when both $M$ and $\Mp$ are hypersurfaces, we also prove the convergence of the associated reflection function when $M$ is assumed to be only minimal. This allows us to derive a new Artin type approximation theorem for formal maps of generic full rank.
\end{abstract}

\maketitle

\section{Introduction}\label{int}

In this paper, we study some properties of formal mappings generically of full rank between generic submanifolds in complex space. A formal holomorphic mapping $H:(\C^N,p)\rightarrow(\C^N,p^\p)$ with $p,p^\p\in\C^N$ is a power series mapping in $z-p$, which satisfy $H(p)=p^\p$. We say that $H$ is generically of full rank if $\J H$, the determinant of the Jacobian matrix of $H$, does not vanish identically (as a formal power series). We define in the usual way the notion of formal mapping between real-analytic generic submanifolds of $\C^N$ (see Section \ref{s:prelim}).

The main problem we were interested in lies in the convergence of formal mappings generically of full rank between real-analytic generic submanifolds. The first result of this type follows from the work of Chern and Moser \cite{CM74} who proved that any formal biholomorphism between real-analytic Levi-nondegenerate hypersurfaces is convergent. The concept of holomorphic nondegeneracy introduced by Stanton \cite{Sta96} turns out to be a crucial condition in understanding the question. Recall that a generic submanifold $M\subset\C^N$ is said to be holomorphically nondegenerate (at $p$) if there exists no non-trivial holomorphic vector field tangent to $M$ (near $p$). The importance of this condition for mapping problems can be seen through the work of Baouendi, Ebenfelt and Rothschild. In the paper \cite{BER97}, the holomorphic nondegeneracy condition is shown to be a necessary condition for the convergence of all formal equivalences between real-analytic generic submanifolds in $\C^N$. In the paper \cite{BMR02}, Baouendi, Mir and Rothschild show that this condition is sufficient to establish the convergence of finite mappings when the source manifold is furthermore assumed to be minimal in the sense of Tumanov (see \cite[\S1.5]{BERbook}).

In this paper, we extend this result to arbitrary formal mappings of generic full rank.

\bt\label{t:cvHgenecase}
Let $(M,p)$, $(\Mp,p^\p)$ be two germs of real-analytic generic manifolds of codimension $d$ in $\C^{n+d}$. Assume that $M$ is minimal at $p$ and that $\Mp$ is \hn \ at $p^\p$. Then any formal holomorphic mapping $H:(\C^{n+d},p)\rightarrow(\C^{n+d},p^\p)$  sending $M$ into $\Mp$ with $\J H\nequiv0$ is convergent.
\et
Our proof of Theorem \ref{t:cvHgenecase} provides also a simpler proof of the analogous result of \cite{BMR02} in the invertible case.
Let us also mention that Theorem \ref{t:cvHgenecase} was proved earlier by Mir \cite{Mir02} in the special case where $\Mp$ is a real-algebraic generic manifold.

As a consequence of Theorem \ref{t:cvHgenecase}, we derive the following convergence result for arbitrary formal maps between real-analytic hypersurfaces.

\bt\label{t:thcourbehol}
Let $(M,p)$, $(\Mp,p^\p)$ be two germs of real-analytic hypersurfaces in $\C^{n+1}$. Assume that $M$ is minimal and \hn \ at $p$ and that $\Mp$ does not contain any holomorphic curve through $p^\p$. Then any formal holomorphic mapping sending $(M,p)$ into $(\Mp,p^\p)$ is convergent. 
\et

Under the more restrictive condition that $M$ is essentially finite at $p$, Theorem \ref{t:thcourbehol} was established by Baouendi, Ebenfelt and Rothschild in \cite{BER00}. In the hypersurface case, we will in fact prove the convergence of the so-called \emph{reflection function} (see \cite{Hua96}) when the source manifold is merely minimal (see Theorem \ref{t:cvR}). This generalizes a result obtained in \cite{Mir00} for formal biholomorphisms and, also allows us to derive the following new Artin type approximation theorem for formal maps of generic full rank.

\bt\label{t:artinhyp}
Let $(M,p)$, $(\Mp,p^\p)$ be two germs of real-analytic \hyp s in $\C^{n+1}$ and $H:(\C^{n+1},p)\rightarrow (\C^{n+1},p^\p)$ be a formal holomorphic mapping sending $M$ into $\Mp$. Assume that $\J H\nequiv 0$ and that $M$ is minimal at $p$. Then for any positive integer $k$, there exists a germ of holomorphic map $H^k:(\C^{n+1},p)\rightarrow (\C^{n+1},p^\p)$ sending $M$ into $\Mp$ which agrees up to order $k$ with $H$ at $p$.
\et

Theorem \ref{t:artinhyp} extends, in the hypersurface case, a similar result of \cite{BMR02} obtained for the more restrictive class of finite mappings. For a more detailed account on recent results on convergence and approximation of formal mappings between submanifolds in complex space, we refer the reader to the survey \cite{MMZ03}.

Our proof of Theorem 1.1 is in essence different from that of Baouendi, Mir and Rothschild for finite mappings \cite{BMR02}. Indeed, in this paper, thanks some mapping identities due to Juhlin \cite{Juh08}, we are able to prove directly the convergence of the formal mapping $H$ (and of its derivatives) along the iterated Segre sets, which strongly contrasts with the methods of \cite{BMR02}. For the hypersurface case, as mentioned earlier, we will prove the convergence of the reflection function for formal maps of generic full rank by combining some arguments of Mir \cite{Mir00} and Juhlin \cite{Juh08}. The Artinian theorem (Theorem \ref{t:artinhyp}) follows directly thanks to a suitable application of Artin's approximation theorem \cite{Artin68} and an argument from \cite{BMR02}.

The paper is organized as follows. In Section \ref{s:prelim}, we recall some basic facts on normal coordinates for generic submanifolds in $\C^N$ and on formal mappings sending real-analytic generic submanifolds into each other. We also prove several preliminary results. Section \ref{s:cvR} is devoted to the proof of the convergence of the reflection function when the source manifold is a minimal hypersurface. Section \ref{s:generalcase} contains the proofs of Theorem \ref{t:cvHgenecase} and Theorem \ref{t:thcourbehol}. To conclude the paper, we prove Theorem \ref{t:artinhyp}.

\section{Preliminaries}\label{s:prelim}

In this section, we first recall the basic definitions we need in this paper, then we establish a convergence lemma which will be used in the next sections.\\

Let $(M,p)$ and $(\Mp,p^\p)$ be germs of real-analytic generic manifolds of codimension $d$ in $\C^{n+d}$. Without loss of generality we may assume that $p=p^\p=0$. It is well known (see for instance \cite[\S4.2]{BERbook}) that there exist normal coordinates defining $M$ and $\Mp$ near 0. It means, in the case of $M$ for example, that there exists $Q\in(\C\{z,\c,\t\})^d$ and $U$, $V$ two open neighborhoods of $0$ in $\C^n$ and $\C^d$ respectively such that $Q$ is convergent on $U\times \bar{U}\times\bar{V}$ and $$M\cap\left(U\times V\right)=\{(z,w)\in U\times V, w=Q(z,\zb,\wb)\}.$$
Moreover, the mapping $Q$ satisfies  
\beq\label{Qnormal}
Q(0,\c,\t)\equiv Q(z,0,\t)\equiv\t
\eeq
and, since $M$ is a real submanifold, we have the following mapping identity
\beq\label{realcond}
Q(z,\c,\Qb(\c,z,\t))\equiv\t.
\eeq
Let us recall that $M$ is said to be minimal at $p$ if there is no germ of a real submanifold $S\subset M$ through $p$ such that the complex tangent space of $M$ at $q$ is tangent to $S$ at every $q\in S$ and $\dim_\mathbb R S<\dim_\mathbb R M$.

We denote by $\mc{M}$ the complexification of $M$ which is the complex submanifold of $\C^{2(n+d)}$ given by 
\beq
\{\left(z,Q\left(z,\c,\t\right),\c,\t\right) , (z,\c,\t)\in U\}\subset\C^{n+d}\times\C^{n+d}.
\eeq
We define in a similar way $\Qp(\zp,\cp,\tp)$ and $\mc{M}^\p$ for $\Mp$. We write the Taylor expansion
\beq\label{Qexpa}
\Qp(\zp,\cp,\tp)=\sum_{\a\in\N^{n}} \Qpa (\cp,\tp){\zp}^\a.
\eeq

Let $H:(\C^n\times\C^d,0)\rightarrow (\C^n\times\C^d,0)$ be a formal holomorphic mapping. We write $H(z,w)=(F(z,w),G(z,w))$ where $F$ and $G$ are respectively with values in $\C^n$ and $\C^d$. 

We say that $H$ sends $M$ into $\Mp$ if the following identity holds at the formal power series level: 
\beq\label{includ}
\Qp(F(z,Q(z,\c,\t)),\Hb(\c,\t))=G(z,Q(z,\c,\t)).
\eeq

We will need the notion of \emph{generic rank} of a formal mapping $h\in(\C[[x]])^q$, $x=(x_1,\ldots,x_r)$. It is the rank of the Jacobian matrix $\deri{x}{h}$ as a matrix with coefficients in the quotient field of $\C[[x]]$. If $h$ is convergent, it coincides with the classical generic rank of a convergent power series mapping.

From now on, we assume that $M$ and $\Mp$ are two real-analytic generic manifolds in $\C^{n+d}$ given in normal coordinates as above. We also assume in what follows that $H=(F,G)$ is a formal mapping sending $M$ into $\Mp$ such that $\J H\nequiv 0$.

We recall a well known fact coming from elementary linear algebra and $\eqref{includ}$.
\bl

In the above situation, the condition $\J H\nequiv 0$ implies  
\beq\label{jac<>0}
\det \dop{z}\left(F\left(z,Q\left(z,\c,\t\right)\right)\right)\nequiv 0.
\eeq

\el

We need the following lemma which can be found in \cite{Juh08} and whose proof is mainly based on a suitable differentiation of the identity $\eqref{includ}$.

\bl\cite[Proposition 5.2]{Juh08}\label{exprQpa}
In the above situation, there exist $\eta_0\in\N^n$ and $\delta_0\in\N^d$ such that, for $p_0=|\eta_0|+|\delta_0|$
\beq
\dbldopz{z}{\eta^\p}{\t}{\delta^\p} \det\dop{z} \left(F(z,Q(z,\c,\t))\right)\equiv 0, \textrm{ for } |\eta^\p|+|\delta^\p|<p_0,
\eeq

\beq
\dbldopz{z}{\eta_0}{\t}{\delta_0} \det\dop{z} \left(F(z,Q(z,\c,\t))\right)\nequiv 0
\eeq
and, for all $\a\in\N^n$, $\b\in\N^d$, such that $|\a|+|\b|> 0$ there exist universal polynomials of their arguments $R_{\a,\b,\eta_0,\delta_0}$ such that, for any $r\in\{1,\ldots,d\}$,

\beq\label{eqQpa}
\left.\mdop{\t}{\b}\right|_{\t=0}\Qpar(\Hb(\c,\t))=\frac{R_{\a,\b,\eta_0,\delta_0}\left(\left(\dbldopz{z}{\eta}{\t}{\delta} H^{(r)}(z,Q(z,\c,\t))\right)_{|\eta|+|\delta|\leq(|\a|+|\b|)(p_0+1)}\right)}{\left(\dbldopz{z}{\eta_0}{\t}{\delta_0} \det \dop{z}\left(F(z,Q(z,\c,\t))\right)\right)^{2(|\a|+|\b|)-1}}
\eeq
where $\Qpar$ is the $r$-th component of $\Qpa$, and $H^{(r)}(z,w)=(F(z,w),G^r(z,w))$ with $G^r$ being the $r$-th component of $G$.
\el

Lemma \ref{exprQpa} is useful to show the following convergence lemma.

\bl\label{l:Qpacv}
In the above situation, for all $\b\in\N^d$, $\a,\g\in\N^n$, $\left.\dbldop{\c}{\g}{\t}{\b}\right|_{\t=0}\Qpa(\Hb(\c,\t))$ is convergent.
\el

\bpf
We fix $\g\in\N^n$. For $|\a|=|\b|=0$, by setting $z=0$, $\t=0$ in equation $\eqref{includ}$, we obtain, thanks to \eqref{Qnormal}, $\Qp(0,\Hb(\c,0))=G(0)$. From the taylor expansion of $\Qp$ given by \eqref{Qexpa}, we get that $\Qp_0(\Hb(\c,0))$ is constant.\\

For $|\a|+|\b|> 0$, we will show that each component of $\left.\dbldop{\c}{\g}{\t}{\b}\right|_{\t=0}\Qpa(\Hb(\c,\t))$ is convergent. For this, we apply $\mdop{\c}{\g}$ to the expression \eqref{eqQpa} of $\left.\mdop{\t}{\b}\right|_{\t=0}\Qpar(\Hb(\c,\t))$ and we obtain that $\left.\dbldop{\c}{\g}{\t}{\b}\right|_{\t=0}\Qpar(\Hb(\c,\t))$ is a ratio of two formal power series. In this ratio, the numerator and the denominator are polynomials in the derivatives of $Q$ which converge and in the derivatives of $H$ at the point $(0,Q(0,\c,0))=0$. So $\left.\dbldop{\c}{\g}{\t}{\b}\right|_{\t=0}\Qpar(\Hb(\c,\t))$ is a ratio of two convergent power series whose denominator does not vanish identically. Since $\left.\dbldop{\c}{\g}{\t}{\b}\right|_{\t=0}\Qpar(\Hb(\c,\t))$ is a formal power series, it is convergent.
\epf

\section{Convergence of the reflection function in the hypersurface case}\label{s:cvR}

In this section, we establish a convergence result for the so-called reflection function of a mapping between minimal real-analytic \hyp s whose Jacobian is not identically zero. We keep the notation of Section \ref{s:prelim} with $d=1$ and still use normal coordinates.

Let $\calR$ be the formal holomorphic mapping called {\em reflection function} (\cite{Hua96}, \cite{Mir00}) and defined by 
\beq\label{defR}
\calR(\zp,\c,\t):=\Qp(\zp,\Hb(\c,\t)),
\eeq
where $(z,\c,\t)\in\C^n\times\C^n\times\C$.
The aim of this section is to establish the following result: 

\bt\label{t:cvR}
Let $M$, $\Mp$ be two real-analytic \hyp s in $\C^{n+1}$ through 0 and $H:(\C^n\times\C,0)\rightarrow (\C^n\times\C,0)$ be a formal holomorphic mapping sending $M$ into $\Mp$ with $\J H\nequiv 0$. If $M$ is minimal at 0 then the associated reflection function $\calR$ given by $\eqref{defR}$ is a convergent power series.
\et

To do this, we follow the different steps of \cite{Mir00}. First, we establish the following result.

\bp\label{p:RcvS1}

Let $M$, $\Mp$ be two real-analytic \hyp s in $\C^{n+1}$ through 0 and $H:(\C^n\times\C,0)\rightarrow (\C^n\times\C,0)$ be a formal holomorphic mapping sending $M$ into $\Mp$. If $\J H\nequiv 0$, then for any $\g\in\N^n$ and $\b\in\N$, the formal holomorphic map 
$$(\zp,\c,\t)\in(\C^n\times\C^n\times\C,0)\mapsto\left.\dbldophyp{\c}{\g}{\t}{\b}\right|_{\t=0}\calR(\zp,\c,\t)$$ is convergent in a neighborhood of 0.
\ep

To prove this proposition, we will use the following convergence lemma whose proof can be found in \cite{Mir00}.

\bl\cite[Lemma 4.1]{Mir00}\label{l:cvborn}
Let $(u_i(t))_{i\in I}$  be a family of convergent power series in $\C\{t\}$, $t=(t_1,\ldots,t_q)$, $q\in\N^*$. Let also $(\mc{K}_i(\z))_{i\in I}$ be a family of convergent power series in $\C\{\z\}$, $\z=(\z_1,\ldots,\z_r)$, $r\in\N^*$. Assume that 
\ben
\item There exists $R>0$ such that the radius of convergence of any $\mc{K}_i$, $i\in I$, is at least R.
\item For all $\z\in\C^r$ with $|\z|<R$, $|\mc{K}_i(\z)|\leq C_i$ with $C_i>0$
\item There exists $V(t)\in(\C[[t]])^r$, $V(0)=0$, such that $\mc{K}_i\circ V(t)=u_i(t)$ for all $i\in I$.
\een
Then, there exists $R^\p>0$ such that the radius of convergence of any $u_i$, $i\in I$, is at least $R^\p$ and such that for all $t\in\C^q$ with $|t|<R^\p$, $|u_i(t)|\leq C_i$.
\el

\bpf[Proof of Proposition {\rm\ref{p:RcvS1}}]

First we will give a proof in the case $\b=|\g|=0$. In order to show that $\calR(\zp,\c,0)$ is convergent, we recall that 

\beq
\calR(\zp,\c,0)=\Qp(\zp,\Hb(\c,0))=\sum_{\a\in\N^{n}} \Qpa (\Hb(\c,0)){\zp}^\a.
\eeq
So, as in the proof of Proposition 5.1 of \cite{Mir00}, we see it suffices to find positive constants $r$ and $R_0$ such that the radius of convergence of the familly $(\Qpa(\Hb(\c,0)))_{\a}$ is at least $r$ and for any $\c$ with norm smaller than $r$ and any $\a\in\N^n$ we have

\beq\label{estimQpa}
|\Qpa(\Hb(\c,0))|\leq R_0^{|\a|+1}.
\eeq
This is possible using the following property which comes from the holomorphy of $\Qp$: there exist $r_1$ and $R$ positive constants such that, for any $\a,\mu\in\N^n$, $\nu\in\N$ and for any $(\cp,\tp)\in\C^{n+1}$ with $|(\cp,\tp)|< r_1$ we have
\beq\label{cauchyestim}
\left|\dblderihyp{\cp}{\mu}{\tp}{\nu}{\Qpa}(\cp,\tp)\right|\leq(\mu,\nu)!R^{|\a|+|\mu|+\nu+1},
\eeq
where $(\mu,\nu)$ is the concatenation of $\mu$ and $\nu$. 

Thus, we can apply Lemma \ref{l:cvborn} to prove $\eqref{estimQpa}$. Indeed, 
\begin{itemize}
\item[-] For any $\a\in\N^n, \Qpa\left(\cp,\tp\right)$ is convergent with radius of convergent at least $r_1$, 
\item[-] by taking $\nu=|\mu|=0$ in $\eqref{cauchyestim}$, we have that for any $(\cp,\tp)\in\C^{n+1}$, with $|(\cp,\tp)|<r_1$, $|\Qpa(\cp,\tp)|\leq R^{|\a|+1}$, 
\item[-] $\Qpa(\Hb(\c,0))$ is convergent thanks to Lemma \ref{l:Qpacv}.
\end{itemize}
Therefore, the proof in the case $\b=|\g|=0$ is complete since $\eqref{estimQpa}$ holds with $R_0=R$.\\

Now, we use the same method to prove the proposition in the general case. We fix $\g\in\N^n$ and $\b\in\N$ with $\b+|\g|>0$. If we write the Taylor expansion of $\calR$ in $\zp$, 
\beq\label{expanR}
\calR(\zp,\c,\t)=\Qp(\zp,\Hb(\c,\t))=\sum_{\a\in\N^n} \Qpa (\Hb(\c,\t)){\zp}^\a,
\eeq
we have
\beq
\trpderihyp{\zp}{\a}{\c}{\g}{\t}{\b}{\calR}(\zp,\c,\t)=\dbldophyp{\c}{\g}{\t}{\b}\left( \mderi{\zp}{\a}{\Qp}\left(\zp,\Hb\left(\c,\t\right)\right) \right).
\eeq
So, to prove the convergence of $\left.\dbldophyp{\c}{\g}{\t}{\b}\right|_{\t=0}\calR(\zp,\c,\t)$, it suffices to prove that the family defined by 
\beq\label{defpsi}
\psi_{\a}^{\b,\g}(\c):=\trpderihyp{\zp}{\a}{\c}{\g}{\t}{\b}{\calR}(0,\c,0)=\a!\left.\dbldophyp{\c}{\g}{\t}{\b}\right|_{\t=0}\Qpa(\Hb(\c,\t))
\eeq
satisfy the following condition: one may find $r_{\b,\g}>0$ and $R_{\b,\g}>0$ such that the radius of convergence of the familly $\psi_{\a}^{\b,\g}$ is at least $r_{\b,\g}$ and such that for any $\a\in\N^n$ and any $\c\in\C^n$ with $|\c|<r_{\b,\g}$, we have the following estimate 
\beq\label{estimpsi}
|\psi_{\a}^{\b,\g}(\c)|\leq \a! R_{\b,\g}^{|\a|+1}.
\eeq
Let's prove $\eqref{estimpsi}$ using Lemma \ref{l:cvborn}. We fix $\a\in\N^n$. First, it comes from $\eqref{defpsi}$ and Lemma \ref{l:Qpacv} that $\psi_{\a}^{\b,\g}(\c)$ is convergent. Then, for any multiindex $\mu\in\N^n$, $\nu\in\N$ with $|\mu|\leq|\g|$, $\nu\leq\b$ there exists a universal polynomial $P_{\nu,\mu}^{\b,\g}$ such that
\beq
\psi_{\a}^{\b,\g}(\c) = 
\a!\sum_{\begin{subarray}{c}|\mu|\leq|\g| \\ \nu\leq\b\end{subarray}}P_{\nu,\mu}^{\b,\g}\left(\left(\dblderihyp{\c}{\delta}{\t}{\eta}{\Hb}(\c,0)\right)_{\begin{subarray}{c}1\leq|\delta|\leq|\g| \\ 1\leq\eta\leq\b\end{subarray}}\right)\dblderihyp{\c}{\mu}{\t}{\nu}{\Qpa}(\Hb(\c,0)).
\eeq
We define a convergent power series mapping of the variables $\left(\left(\L_{\delta,\eta}\right)_{\begin{subarray}{c}1\leq|\delta|\leq|\g| \\ 1\leq\eta\leq\b\end{subarray}},\cp,\tp\right)\in\C^{N_{\b,\g}}\times\C^{n+1}$ with $N_{\b,\g}=(n+1){\rm Card}\{(\delta,\eta)\in\N^n\times\N,1\leq|\delta|\leq|\g|,1\leq\eta\leq\b\}$, as follows 
\beq
h_{\a}^{\b,\g}\left(\left(\L_{\delta,\eta}\right)_{\begin{subarray}{c}1\leq|\delta|\leq|\g| \\ 1\leq\eta\leq\b\end{subarray}},\cp,\tp\right):=
\a!\sum_{\begin{subarray}{c}|\mu|\leq|\g| \\ \nu\leq\b\end{subarray}}P_{\nu,\mu}^{\b,\g}\left(\left(\L_{\delta,\eta}+\dblderihyp{\c}{\delta}{\t}{\eta}{\Hb}(0)\right)_{\begin{subarray}{c}1\leq|\delta|\leq|\g| \\ 1\leq\eta\leq\b\end{subarray}} \right)\dblderihyp{\c}{\mu}{\t}{\nu}{\Qpa}(\cp,\tp).
\eeq
From this and using $\eqref{cauchyestim}$, we observe that there exists a constant $R_{\b,\g}$ such that, for $\left(\L_{\delta,\eta}\right)_{\begin{subarray}{c}1\leq|\delta|\leq|\g| \\ 1\leq\eta\leq\b\end{subarray}}$ with norm smaller than 1 in $\C^{N_{\b,\g}}$ and $(\cp,\tp)$ with norm smaller than $r$ in $\C^{n+d}$, we have, for any $\a\in\N^n$
\beq
\left|h_{\a}^{\b,\g}\left(\left(\L_{\delta,\eta}\right)_{\begin{subarray}{c}1\leq|\delta|\leq|\g| \\ 1\leq\eta\leq\b\end{subarray}},\cp,\tp\right)\right|\leq \a!R_{\b,\g}^{|\a|+1}.
\eeq
So, using the fact that
\beq
h_{\a}^{\b,\g}\left(\left(\dblderihyp{\c}{\delta}{\t}{\eta}{\Hb}(\c,0)-\dblderihyp{\c}{\delta}{\t}{\eta}{\Hb}(0)\right)_{\begin{subarray}{c}1\leq|\delta|\leq|\g| \\ 1\leq\eta\leq\b\end{subarray}},\Hb(\c,0)\right)=\psi_{\a}^{\b,\g}(\c)
\eeq
and Lemma \ref{l:cvborn}, we obtain the estimates $\eqref{estimpsi}$. Thus the proof of Proposition \ref{p:RcvS1} is complete.
\epf

We now need the two following lemma and proposition. The proof of the lemma, which is a consequence of Artin's approximation theorem \cite{Artin68}, can be found in \cite{Mir00}. The proposition is due to Gabrielov \cite{Gab73} and proved in e.g. \cite{BM88}.

\bl\cite[Lemma 6.1]{Mir00}\label{l:artinderi}
Let $\mc{T}(x,u)\in(\C[[x,u]])^r$, $x\in\C^q$, $u\in\C^s$with $\mc{T}(0)=0$. Assume that $\mc{T}(x,u)$ satisfies an identity in the ring $\C[[x,u,y]]$, $y\in\C^q$, of the form 
$$\phi(\mc{T}(x,u);x,u,y)=0$$
where $\phi\in\C[[W,x,u,y]]$ with $W\in\C^r$. Assume, furthermore, that for any multiindex $\b\in\N^q$, the formal power series $\left[\mderi{y}{\b}{\phi}\left(W;x,u,y\right)\right]_{y=x}$ is convergent. Then, for any given positive integer $k$, there exists an $r$-uple of convergent power series $\mc{T}^k(x,u)$ such that $\phi(\mc{T}^k(x,u);x,y,u)=0$ in $\C[[x,u,y]]$ and such $\mc{T}^k(x,u)$ agrees up to order $k$ at 0 with $\mc{T}(x,u)$.
\el

\bp\label{p:cv+genesub}
Let $\mc{J}(x)\in(\C\{x\})^r$, $x\in\C^k$, $k,r\geq 1$ $\mc{J}(0)=0$ and $\mc{V}(t)\in\C[[t]]$, $t\in\C^r$. If $\mc{V}\circ\mc{J}$ is convergent and the generic rank of $\mc{J}$ is $r$, then $\mc{V}$ is convergent.
\ep

We will also need the definition and properties of the \emph{Segre set mappings} for the hypersurface $M$, only the three first for this section: 
\begin{eqnarray}
&v_1:&z\in(\C^n,0)\mapsto(z,0), \\ \nonumber
&v_2:&(z,\xi)\in(\C^{2n},0)\mapsto(z,Q(z,\xi,0)), \\ \nonumber
&v_3:&(z,\xi,\z)\in(\C^{3n},0)\mapsto(z,Q(z,\xi,\Qb(\xi,\z,0))). \nonumber
\end{eqnarray}
The minimality condition of $M$ in Theorem \ref{t:cvR} allows us to say that the generic rank of $v_2$ (and $v_3$) is $n+1$.

\bpf[Proof of Theorem {\rm\ref{t:cvR}}]
By setting $\t=\Qb(\c,z,Q(z,\xi,0))$ in $\eqref{includ}$ and using $\eqref{realcond}$, we obtain 
\beq
\Qp\left(F\circ v_2(z,\xi),\Hb\circ \bar{v}_3(\c,z,\xi)\right)=G\circ v_2(z,\xi)
\eeq
This means that $\calR$ satisfy the following equation  
\beq\label{Req}
\calR(F\circ v_2(z,\xi),\bar{v}_3(\c,z,\xi))-G\circ v_2(z,\xi)\equiv0.
\eeq
We want to apply Lemma \ref{l:artinderi} to the equation $\eqref{Req}$ with $x=\xi$, $u=z$, $y=\c$, $\mc{T}(x,u)=H\circ v_2(\c,\xi)$, $W=(\l,\mu)\in\C^n\times\C$ and 
$$\phi((\l,\mu);\xi,z,\c)=\calR(\l,\bar{v}_3(\c,z,\xi))-\mu.$$
For this, we have to verify the convergence of all the derivatives of $\phi$ with respect to $\c$ evaluated in $\c=\xi$. But, those derivatives involve only derivatives of $v_3$ which are convergent and expressions of the form $\dblderi{\c}{\g}{\t}{\b}{\calR}(\l,\bar{v}_3(\c,z,\xi))|_{\c=\xi}$.
As, from the mapping identity $\eqref{realcond}$, $v_3(\xi,z,\xi)=v_1(\xi)=(\xi,0)$, Proposition \ref{p:RcvS1} allows us to apply Lemma \ref{l:artinderi}. Thus, for any positive integer $k$ there exists $\mc{T}^k(z,\xi)=({\mc{T}^\p}^k(z,\xi),\mc{T}^k_n(z,\xi))$ convergent power series agreeing up to order $k$ at 0 with $H\circ v_2(z,\xi)$ and such that 
$$\calR\left({\mc{T}^\p}^k(z,\xi),\bar{v}_3(\c,z,\xi)\right)\equiv\mc{T}^k_n(z,\xi)$$
To show that $\calR$ is convergent it suffices, using Proposition \ref{p:cv+genesub}, to show that, for $k$ well chosen, the map $(z,\c,\xi)\in\C^{3n}\mapsto\left({\mc{T}^\p}^k(z,\xi),\bar{v}_3(\c,z,\xi)\right)\in\C^{2n+1}$ is of generic rank $2n+1$.

From the two following facts: 

\ben
\item $(z,\c,\xi)\in\C^{3n}\mapsto(v_2(z,\xi),\bar{v}_3(\c,z,\xi))\in\mc{M}$ is of generic rank $2n+1$ since $M$ is minimal,
\item $(z,w,\c,\t)\in\mc{M}\mapsto(F(z,w),\c,\t)$ is of generic rank $2n+1$ because $F$ is of generic full rank $n$;
\een
we deduce that $(z,\c,\xi)\mapsto(F(v_2(z,\xi)),\bar{v}_3(\c,z,\xi))$ is of generic rank $2n+1$ in $(\C[[z,\c,\xi]])^{2n+1}$. But $\mc{T}^k(z,\xi)$ agrees up to order $k$ with $H\circ v_2(z,\xi)$, so for $k$ large enough the generic rank of $\left({\mc{T}^\p}^k(z,\xi),\bar{v}_3(\c,z,\xi)\right)$ is $2n+1$. This completes the proof of Theorem \ref{t:cvR}.
\epf

\section{The case of generic manifolds of arbitrary codimension}\label{s:generalcase}

The aim of this section is to show Theorem \ref{t:cvHgenecase}. The arguments used to obtain this result are inspired by an analogous result in the finite jet determination problem that can be found in \cite{Juh08}.

\subsection{A convergence lemma}

In this paragraph, we will state and prove a lemma inspired from \cite{Juh08}. 

\bl\label{l:lemdecv}
Let $P(x,Y)\in(\C\{x,Y\})^N$ with $x\in\C^{n_1}$, $Y\in\C^N$, $\phi_0(x,t)\in(\C[[x,t]])^N$, $\phi(0)=0$, with $t\in\C^d$ such that 
$$\det\deri{Y}{P}(x,\phi_0(x,t))\nequiv 0.$$
Assume that for every $\b\in\N^d$ $\mdopz{t}{\b} P(x,\phi_0(x,t))$ is a convergent power series, then for every $\b\in\N^d$ $\mderi{t}{\b}{\phi_0}(x,0)$ is a convergent power series too.
\el

To prove this result we need the following proposition from \cite{Juh08}.

\bp\cite[Proposition 3.1]{Juh08},\label{p:propdejuh}
Let $P(x,Y)\in(\C[[x,Y]])^N$ with $x\in\C^{n_1}$, $Y\in\C^N$, $\phi_0(x,t)\in(\C[[x,t]])^N$, $\phi_0(0)=0$, with $t\in\C^d$ such that 
\beq\label{detP<>0}
\det\deri{Y}{P}(x,\phi_0(x,t))\nequiv 0.
\eeq
Let $\a_0\in\N^{n_1}$, $\b_0\in\N^d$ such that 

\beq
\dbldopz{x}{\a_0}{t}{\b_0}\det\deri{Y}{P}(x,\phi_0(x,t))\neq0.
\eeq
Then, for any power series $\phi(x,t)$ satisfying 
\ben
\item for every $(\a,\b)\in\N^{n_1}\times\N^d$ with $|\a|\leq|\a_0|$ and $|\b|\leq|\b_0|$:  $\dblderi{x}{\a}{t}{\b}{\phi}\left(0\right)=\dblderi{x}{\a}{t}{\b}{\phi_0}\left(0\right)$,
\item for some $k\in\N$ we have $\mdopz{t}{\b}P(x,\phi(x,t))=\mdopz{t}{\b}P(x,\phi_0(x,t))$, for $|\b|\leq|\b_0|+k$,
\een

\noindent we have that, for every $\b\in\N^d$ with $|\b|\leq k$,
\beq
\mderi{t}{\b}{\phi}(x,0)=\mderi{t}{\b}{\phi_0}(x,0).
\eeq

\ep

\bpf[Proof of Lemma {\rm \ref{l:lemdecv}}]
For fixed $\b\in\N^d$, we note $P_\b(x):=\mdopz{t}{\b} P(x,\phi_0(x,t))\in(\C\{x\})^N$. From the chain rule, there exists a universal matrix with polynomial coefficients $\left(A_{\nu,\b}\left(\left(\L_\mu\right)_{1\leq|\mu|\leq |\b|}\right)\right)_{|\nu|\leq|\b|}$ where $(\L_\mu)_{1\leq|\mu|\leq |\b|}\in\C^{N_\b}$ and $N_\b=(n+d){\rm Card}\{\mu\in\N^d, 1\leq|\mu|\leq |\b| \}$ satisfying: 
\beq\label{exprederiP}
P_{\b}(x)=\sum_{|\nu|\leq|\b|} \mderi{Y}{\nu}{P}\left(x,\phi_0(x,0)\right)\cdot A_{\nu,\b}\left(\left(\mderi{t}{\mu}{\phi_0}(x,0)\right)_{1\leq|\mu|\leq |\b|}\right).
\eeq
We define convergent power series by setting 
\beq
h_\b((\L_\mu)_{0\leq|\mu|\leq|\b|},x):=\sum_{|\nu|\leq|\b|}\mderi{Y}{\nu}{P}\left(x,\L_0\right) \cdot A_{\nu,\b}\left(\left(\L_\mu\right)_{1\leq|\mu|\leq |\b|}\right)
\eeq
and 
\beq
R_\b((\L_\mu)_{0\leq|\mu|\leq|\b|},x):=h_\b((\L_\mu)_{0\leq|\mu|\leq|\b|},x)-P_\b(x).
\eeq
Now we fix $L\in\N$ and we consider 
\beq
R\left(\left(\L_\mu\right)_{0\leq|\mu|\leq L},x\right):=\left(\left( R_\beta(\L_{0\leq|\mu|\leq|\b|}),x\right)\right)_{0\leq|\b|\leq L}.
\eeq
By definition of $R$, it is a convergent power series and it satisfies
\beq
R\left(\left(\mderi{t}{\mu}{\phi_0}(x,0)\right)_{0\leq|\mu|\leq L},x\right)\equiv 0.
\eeq
So, by Artin's approximation theorem \cite{Artin68}, for any positive integer $K$, there exist a convergent mapping  $\left(\phi_K^{\mu, L}(x)\right)_{0\leq|\mu|\leq L}\in(\C\{x\})^{N_L}$ where $N_L=(n+d){\rm Card}\{\mu\in\N^d, |\mu|\leq L\}$, which agrees up to order $K$ with $\left(\mderi{t}{\mu}{\phi_0}(x,0)\right)_{0\leq|\mu|\leq L}$  at 0 and such that 
\beq\label{equationphiKmuL}
R\left(\left(\phi_K^{\mu, L}(x)\right)_{0\leq|\mu|\leq L},x\right)=0.
\eeq
We define the following convergent power series $\phi(x,t):=\sum_{|\mu|\leq L} \frac{\phi_K^{\mu, L}(x)}{\mu!} t^\mu$.
By definition we have that, for any $\mu\in\N^d$ with $|\mu|\leq L$, $\mderi{t}{\mu}{\phi}(x,0)=\phi_K^{\mu, L}(x)$. Thus, we obtain that, for any $(\a,\mu)\in\N^{n_1}\times\N^d$ with $|\a|\leq K$ and $|\mu|\leq L$, 
\beq\label{deriphi}
\dblderi{x}{\a}{t}{\mu}{\phi}(0)=\dblderi{x}{\a}{t}{\mu}{\phi_0}(0).
\eeq
So, from $\eqref{equationphiKmuL}$, $R\left(\left(\mderi{t}{\mu}{\phi}(x,0)\right)_{0\leq|\mu|\leq L},x\right)=0$, which is equivalent to 

\beq\label{Pdephi=Pdephi_0}
\mdopz{t}{\b} P(x,\phi(x,t))=P_\b(x)=\mdopz{t}{\b} P(x,\phi_0(x,t)),
\eeq
for any $\b\in\N^d$ with $|\b|\leq L$.

Now, we want to apply Proposition \ref{p:propdejuh} by making a good choice of $K$ and $L$. If $\a_0$ and $\b_0$ are like in  Proposition \ref{p:propdejuh} (this is possible because of $\det\deri{Y}{P}\left(x,\phi_0\left(x,t\right)\right)\nequiv0$), we take $K=|\a_0|$ and $L=|\b_0|+k$ for any fixed given $k>0$. From $\eqref{deriphi}$ and $\eqref{Pdephi=Pdephi_0}$, the formal mappings $\phi_0$ and $\phi$ satisfy the hypothesis of Proposition \ref{p:propdejuh} so that we have that, for any $\b\in\N^d$ with $|\b|\leq k$,
\beq
\mderi{t}{\b}{\phi_0}(x,0)=\mderi{t}{\b}{\phi}(x,0).
\eeq
As $\phi$ is a convergent power series, $\mderi{t}{\b}{\phi_0}(x,0)$ is convergent for $|\b|\leq k$. This is true for any choice of $k$, so the proof of Lemma \ref{l:lemdecv} is complete.
\epf

\subsection{Proof of Theorem \ref{t:cvHgenecase}}

Let us recall the setting of Theorem \ref{t:cvHgenecase}. We let $(M,p)$ and $(\Mp,p^\p)$ be two germs of real-analytic generic manifolds of codimension $d\geq1$ in $\C^ {n+d}$ with $M$ minimal at $p$ and $\Mp$ holomorphically nondegenerate at $p^\p$. We consider $H:(\C^{n+d},p)\rightarrow(\C^{n+d},p^\p)$ a formal mapping which sends $M$ into $\Mp$ with $\J H\nequiv0$. We use the notation of Section \ref{s:prelim} and work in normal coordinates.

To obtain the convergence of $H$, we will establish its convergence along the \emph{Segre sets} (see \cite[Chapter X]{BERbook}). The following result provides the convergence along the first \emph{Segre set}.

\bp\label{p:convHrdS1}
Under the assumptions of Theorem {\rm\ref{t:cvHgenecase}} and with the above notation, for any $\b\in\N^d$, $\mderi{w}{\b}{H}(z,0)$ is convergent.
\ep

\bpf
Since $\Mp$ is \hn, using Stanton's criterion (see \cite{Sta96} or \cite[\S11.3]{BERbook}), it is possible to choose $\a^i\in\N^n$ and $1\leq r_i\leq d$, for $i\in\{1,\dots,n\}$ such that $\det\left(\deri{\cp_l}{\Qp^{r_k}_{\a^k}}(\cp,\tp)\right)_{k,l}\nequiv 0$. Let $P(\cp,\tp)=\left(P^1(\cp,\tp),\ldots,P^{n+d}(\cp,\tp)\right)$ be the convergent power series mapping defined by 

\begin{displaymath}
P^k(\cp,\tp):=\left\{
\begin{array}{ll}
\Qp^{r_k}_{\a^k}(\cp,\tp),& \textrm{ for } 1\leq k\leq n,  \\ \nonumber
\Qp^{k-n}_{0}(\cp,\tp),& \textrm{ for } n+1\leq k\leq n+d. 
\end{array}\right.
\end{displaymath}
Since $\deri{\tp}{Q^\p_0}(\cp,\tp)=I_d$ and $\deri{\cp}{Q^\p_0}(\cp,\tp)=0$, we have $$\l(\cp,\tp):=\det\deri{(\cp,\tp)}{P}(\cp,\tp)=\det\left(\deri{\cp_l}{\Qp^{r_k}_{\a^k}}(\cp,\tp)\right)_{k,l}\nequiv 0.$$
Since $H$ is of generic full rank, $\l(\Hb(\c,\t))\nequiv0$. Moreover, thanks to Lemma \ref{l:Qpacv}, $\mdopz{\t}{\b} P(\Hb(\c,\t))$ is convergent for any $\b\in\N^d$. It follows from Lemma \ref{l:lemdecv} that, for any $\b\in\N^d$, $\mdopz{\t}{\b}\Hb(\c,\t)$ is convergent which is equivalent to $\mderi{w}{\b}{H}(z,0)$ convergent.
\epf

The convergence of $H$ along the \emph{Segre sets} of higher order will be proved using the \emph{Segre set mappings} (see \cite[\S10.4]{BERbook}).

For $j\in\N^*$, we define the convergent power series $U_j:(\C^{nj}\times\C^d,0)\rightarrow (\C^d,0)$ by induction on $j$ as follows:

\beq
U_1\left (z^1;t\right ):=t,
\eeq

\beq
U_{2j}\left (z^1,\c^1,z^2,\c^2,\ldots,z^j,\c^j;t\right ):=U_{2j-1}\left (z^1,\c^1,z^2,\c^2,\ldots,z^j;Q\left (z^j,\c^j,t\right )\right ), 
\eeq

\beq
U_{2j+1}\left (z^1,\c^1,z^2,\c^2,\ldots,z^j,\c^j,z^{j+1};t\right ):=U_{2j}\left (z^1,\c^1,z^2,\c^2,\ldots,z^j,\c^j;\Qb\left (\c^j,z^{j+1},t\right )\right ), 
\eeq
where $z^k,\c^k\in\C^n$ and $t\in\C^d$.
Let $V_j:(\C^{nj}\times\C^d,0)\rightarrow(\C^n\times\C^d,0)$ be the convergent power series mapping defined by
\beq\label{Sj}
V_j\left (z^1,\c^1,z^2,\c^2,\ldots;t\right ):=\left (z^1,U_j\left (z^1,\c^1,z^2,\c^2,\ldots;t\right )\right ).
\eeq
From the mapping identity $\eqref{realcond}$ we deduce, for any $j\geq1$,
\beq\label{realcondSj}
\left(z,Q\left(z^1,\c^1,\Ub_{j+1}\left(\c^1,z^1,\c^2,z^2,\ldots;t\right)\right)\right)=V_j(z^1,\c^2,z^2,\ldots;t).
\eeq
By setting $\t=\Ub_{j+1}\left(\c^1,z^1,\c^2,z^2,\ldots;t\right)$ in $\eqref{includ}$ and using $\eqref{realcondSj}$, we obtain the following power series identity, for any $j\geq 1$,
\beq\label{includsegre}
\Qp\left (F(V_j(z^1,\c^2,z^2,\ldots;t)),\Hb(\Vb_{j+1}(\c^1,z^1,\c^2,z^2,\ldots;t)))=G(V_j(z^1,\c^2,z^2,\ldots;t)\right ).
\eeq
We will prove the following convergence result: 

\bp\label{p:convHrdSj}
Under the assumptions of Theorem {\rm\ref{t:cvHgenecase}} and in the above setting, we have that, for every $\b\in N^d$ and every integer $j\geq 1$ 
$$\mdopz{t}{\b}\Big(H\rd V_j(\c^1,z^1,\c^2,z^2,\ldots;t)\Big) \textrm{ converges. }$$
\ep

We first need the following result whose proof is inspired from \cite{Juh08}.

\bl\label{l:lemcvh}
Let $h(Z)$ be a power series with $Z=(z,w)\in\C^n\times\C^d$ and $j\geq1$ a fixed integer. If for every $\b\in\N^d$ $$\mdopz{t}{\b} \Big(h\rd V_j(z^1,\c^1,\ldots;t)\Big) \textrm{ converges},$$ then for every $\g\in\N^{n+d}$ and every $\b\in\N^d$, $\mdopz{t}{\b}\Big( \left(\mderi{Z}{\g}{h}\right)\rd V_j(z^1,\c^1,\ldots;t)\Big)$ converges.
\el

\bpf
We prove this lemma by induction on $|\g|$. The case $\g=0$ comes from the hypothesis. Now assume the result for any $\g\in\N^{n+d}$, $|\g|=l$ and $\b\in\N^d$, we show the analogous conclusion for $\g^\p\in\N^{n+d}$, $|\g^\p|=l+1$ and any $\b\in\N^d$. We denote 
\beq\label{defRj}
R^j_\g(z^1,\c^1,\ldots;t):=\left(\mderi{Z}{\g}{h}\rd V_j\right)(z^1,\c^1,\ldots;t).
\eeq
The induction hypothesis implies that $\mdopz{t}{\b}R^j_\g$ are convergent. Differentiating $\eqref{defRj}$ with respect to $t$, we obtain
\beq\label{deri t Rj}
\deri{t}{R^j_\g}(z^1,\c^1,\ldots;t)=\left(\frac{\d^{|\g|+1} h}{\d Z^\g \d w}\rd V_j\right)(z^1,\c^1,\ldots;t)\cdot\deri{t}{U_j}(z^1,\c^1,\ldots;t).
\eeq
Since $\deri{\t}{Q}(0)=I_d$, we have $\deri{t}{U_j}(0)=I_d$. And, multiplying $\eqref{deri t Rj}$, by the classical adjoint matrix of $\deri{t}{U_j}(z^1,\c^1,\ldots;t)$, i.e. the transpose of the cofactor matrix denoted by $\adj\left(\deri{t}{U_j}(z^1,\c^1,\ldots;t)\right)$, we obtain 
\beq
\left(\frac{\d^{|\g|+1} h}{\d Z^\g \d w}\rd V_j\right)(z^1,\c^1,\ldots;t)=\frac{\deri{t}{R^j_\g}(z^1,\c^1,\ldots;t)\cdot\adj\left(\deri{t}{U_j}(z^1,\c^1,\ldots;t)\right)}{\det\deri{t}{U_j}(z^1,\c^1,\ldots;t)}.
\eeq
Now, if we apply $\mdopz{t}{\b}$ to the previous identity we get that $\mdopz{t}{\b}\left(\frac{\d^{|\g|+1} h}{\d Z^\g \d w}\rd V_j\right)(z^1,\c^1,\ldots;t)$ is a ratio of two convergent power series whose denominator does not vanish at 0. Therefore $\mdopz{t}{\b}\left(\frac{\d^{|\g|+1} h}{\d Z^\g \d w}\rd V_j\right)(z^1,\c^1,\ldots;t)$ is convergent.
To obtain the other derivatives, we differentiate $\eqref{defRj}$ with respect to $z^1$, we get 
\beq
\left(\frac{\d^{|\g|+1 } h}{\d Z^\g \d z}\rd V_j\right)(z^1,\c^1,\ldots;t)+\left(\frac{\d^{|\g|+1} h}{\d Z^\g \d w}\rd V_j\right)(z^1,\c^1,\ldots;t)\cdot\deri{z^1}{U_j}(z^1,\c^1,\ldots;t)=\deri{z^1}{R^j_\g}(z^1,\c^1,\ldots;t).
\eeq
All the derivatives of $\mdopz{t}{\b}\deri{z^1}{R^j_\g}(z^1,\c^1,\ldots;t)$ are convergent. Indeed, $\mdopz{t}{\b}R^j_\g(z^1,\c^1,\ldots;t)$ is convergent by hypothesis, thus $\dop{z^1}\left(\mdopz{t}{\b}R^j_\g(z^1,\c^1,\ldots;t)\right)$ is convergent too. This implies that $\mdopz{t}{\b}\left(\deri{z^1}{R^j_\g}(z^1,\c^1,\ldots;t)\right)$ is convergent. 

Thus, using the fact that $\mdopz{t}{\b}\left(\frac{\d^{|\g|+1} h}{\d Z^\g \d w}\rd V_j\right)(z^1,\c^1,\ldots;t)$ is convergent by the previously established result, we obtain the convergence of $\mdopz{t}{\b}\left(\frac{\d^{|\g|+1} h}{\d Z^\g \d z}\rd V_j\right)(z^1,\c^1,\ldots;t)$. Consequently, the proof of Lemma \ref{l:lemcvh} is complete. 
\epf

To prove Proposition \ref{p:convHrdSj}, we will also need the following mapping identity result whose proof can be found in \cite{Juh08}.

\bl\cite[Proposition 6.4]{Juh08}\label{l:lem6.4deJuh}
In the above situation and for $j\geq1$ a fixed integer, as $\J H\nequiv0$ we may choose $\g^j\in\N^d$ such that 
\beq
\mdopz{t}{\delta^\p} \det\left(\left.\left(\dop{z^1}\left( F(z^1,Q(z^1,\c^1,\t))\right)\right)\right|_{\t=\Ub_{j+1}(\c^1,z^1,\c^2,\ldots;t)}\right)\equiv 0, \textrm{ for } |\delta^\p|<|\g^j|,
\eeq

\beq
\mdopz{t}{\g^j} \det\left(\left.\left(\dop{z^1}\left( F(z^1,Q(z^1,\c^1,\t))\right)\right)\right|_{\t=\Ub_{j+1}(\c^1,z^1,\c^2,\ldots;t)}\right)\nequiv 0,
\eeq

and, so that for every positive integer $k$, we have
\begin{eqnarray}
\frac{1}{(k\g^j)!}\left.\frac{\d^{k|\g^j|}}{\d t^{k\g^j}}\right|_{t=0}\left\{\det\left(\left.\left( \dop{z^1}\left(F(z^1,Q(z^1,\c^1,\t))\right)\right)\right|_{\t=\Ub_{j+1}(\c^1,z^1,\c^2,\ldots;t)}\right)^k\right\} \\ \nonumber
=\left(\frac{1}{\g^j!}\left.\frac{\d^{|\g^j|}}{\d t^{\g^j}}\right|_{t=0}\det\left(\left.\left( \dop{z^1}\left(F(z^1,Q(z^1,\c^1,\t))\right)\right)\right|_{\t=\Ub_{j+1}(\c^1,z^1,\c^2,\ldots;t)}\right)\right)^k\nonumber.
\end{eqnarray}

Then, for any $\a\in\N^n$ and $\b\in\N^d$ with $|\a|+|\b|>0$ there is a universal polynomial of its arguments $T_{\a,\b,\g^j}$ such that the following holds: for any $j\geq1$ and every $r\in\{1,\ldots,r\}$ , we have the mapping identity 
\beq\label{exprQplem6.4}
\begin{split}
&\mdopz{t}{\b}\Qp^r_{\zp^\a}\left (F\rd V_j(z^1,\c^2,z^2,\ldots;0),\Hb\rd \Vb_{j+1}(\c^1,z^1,\c^2,\ldots;t)\right ) \\ 
&=\frac{T_{\a,\b,\g^j}\left(\left(  \mdopz{t}{\delta}\left(\left.\left(\mdop{z^1}{\eta} H^{(r)}(z^1,Q(z^1,\c^1,\t))\right)\right|_{\t=\Ub_{j+1}(\c^1,z^1,\c^2,\ldots;t)} \right)\right)_{|\eta|+|\delta|\leq (|\a|+|\b|)(2|\g^j|+1)-|\g^j|}\right)}{\left(\frac{1}{\g^j!}\left.\frac{\d^{|\g^j|}}{\d t^{\g^j}}\right|_{t=0}\det\left(\left.\left( \dop{z^1}\left(F(z^1,Q(z^1,\c^1,\t))\right)\right)\right|_{\t=\Ub_{j+1}(\c^1,z^1,\c^2,\ldots;t)}\right)\right)^{2(|\a|+|\b|)-1} }.
\end{split}
\eeq
\el

\bpf[Proof of Proposition {\rm\ref{p:convHrdSj}}]
We prove this proposition by induction on $j$. For $j=1$, the result comes from Proposition \ref{p:convHrdS1}. We assume the proposition is true for fixed $j$ and we prove it for $j+1$. To do this we need the following lemma:

\bl\label{l:convQFSjHSj+1}
For any $\b\in\N^d$, $\a\in\N^n$ and $r\in\{1,\ldots,d\}$, $$\mdopz{t}{\b}\Qp^r_{\zp^\a}\left (F\rd V_j(z^1,\c^2,z^2,\ldots;0),\Hb\rd \Vb_{j+1}(\c^1,z^1,\c^2,\ldots;t)\right )\textrm{ is convergent.}$$
\el

\bpf
We begin by the case $|\a|=|\b|=0$. By setting $t=0$ in $\eqref{includsegre}$, we obtain 
$$\Qp\left (F\rd V_j(z^1,\c^2,z^2,\ldots;0)),\Hb\rd \Vb_{j+1}(\c^1,z^1,\c^2,\ldots;0)\right )=G\rd V_j(z^1,\c^2,z^2,\ldots;0)).$$
The left hand side of the previous equation is exactly the expression whose convergence has to be proven and the right hand side is convergent by the induction hypothesis.

For the case $|\a|+|\b|>0$, from Lemma \ref{l:lem6.4deJuh} we have that

\begin{equation*}
\begin{split} 
&\mdopz{t}{\b}\Qp^r_{\zp^\a}\left (F\rd V_j(z^1,\c^2,z^2,\ldots;0),\Hb\rd \Vb_{j+1}(\c^1,z^1,\c^2,\ldots;t)\right ) \\
&=\frac{T_{\a,\b,\g}\left(\left(  \mdopz{t}{\delta}\left(\left.\left(\mdop{z^1}{\eta} H^{(r)}(z^1,Q(z^1,\c^1,\t))\right)\right|_{\t=\Ub_{j+1}(\c^1,z^1,\c^2,\ldots;t)} \right)\right)_{|\eta|+|\delta|\leq (|\a|+|\b|)(2|\g|+1)-|\g|}\right)}{\left(\frac{1}{\g!}\left.\frac{\d^{|\g|}}{\d t^{\g}}\right|_{t=0}\det\left(\left.\left( \dop{z^1}\left(F(z^1,Q(z^1,\c^1,\t))\right)\right)\right|_{\t=\Ub_{j+1}(\c^1,z^1,\c^2,\ldots;t)}\right)\right)^{2(|\a|+|\b|)-1} }.
\end{split}
\end{equation*}
This is a ratio of two convergent power series. Indeed, the numerator is polynomial in the derivatives of $Q$ which are convergent and expressions of the form  
$$\mdopz{t}{\delta}\dblderi{z}{\eta}{w}{\mu}{H}\left (z^1,Q(z^1,\c^1,\Ub_{j+1}(\c^1,z^1,\c^2,\ldots;t)\right ),$$
i.e. using $\eqref{realcondSj}$, of the form 
$$\mdopz{t}{\delta}\left(\dblderi{z}{\eta}{w}{\mu}{H}\rd V_j\right)(z^1,\c^2,z^2,\ldots;t)$$ 
which are also convergent thanks to the induction hypothesis and Lemma \ref{l:lemcvh}.
\epf

We come back to the proof of Proposition \ref{p:convHrdSj}. As in the proof of Proposition \ref{p:convHrdS1}, from the holomorphic nondegeneracy assumption on $\Mp$,  we may choose $\a^i\in\N^n$ and $1\leq r_i\leq d$ for $i\in\{1,\ldots,n\}$ such that 
\beq\label{hypohn}
\det\left(\deri{\cp_l}{\Qp^{r_k}_{\a^k}}(\cp,\tp)\right)_{k,l}\nequiv 0.
\eeq
We define a power series mapping $P_{j+1}:(\C^{(j+1)n}\times\C^{n+d},0)\rightarrow (\C^{n+d},0)$ in the following way: 
$$
P_{j+1}(\c^1,z^1,\c^2,z^2,\ldots,Y):=\left(
\begin{array}{c}
\Qp^{r_1}_{\zp^{\a^1}}(F(V_j(z^1,\c^2,z^2,\ldots;0)),Y)\\
\vdots\\
\Qp^{r_n}_{\zp^{\a^n}}(F(V_j(z^1,\c^2,z^2,\ldots;0)),Y)\\
\Qp^1(F(V_j(z^1,\c^2,z^2,\ldots;0)),Y)\\
\vdots\\
\Qp^d(F(V_j(z^1,\c^2,z^2,\ldots;0)),Y)
\end{array}\right).
$$
By the induction hypothesis, $P_{j+1}$ is a convergent power series mapping. Moreover, we observe that 
\beq\label{deriPj+1}
\det\deri{Y}{P_{j+1}}\left (\c^1,z^1,\c^2,\ldots;\Hb\rd\Vb_{j+1}(\c^1,z^1,\c^2,\ldots;t)\right )\nequiv0.
\eeq
Indeed, if we set, in $\eqref{deriPj+1}$, $z^l=0$ for $l\geq 1$ and $\c^l=0$ for $l\geq 2$, we obtain, since $\Vb_{j+1}(\c^1,0,\ldots;t)=(\c^1,t)$, $\deri{\cp}{\Qp}(0,\cp,\tp)=0$ and $\deri{\tp}{\Qp}(0,\cp,\tp)=I_d$: 
\beq
\det\deri{Y}{P_{j+1}}\left (\c^1,0,\ldots;\Hb(\c^1,t)\right )=\det\left(\deri{\cp_l}{\Qp_{\zp^{\a^k}}^{r_k}}\left (0,\Hb(\c^1,t)\right )\right)_{k,l}=\det\left(\deri{\cp_l}{\Qp_{\a^k}^{r_k}}\left (\Hb(\c^1,t)\right )\right)_{k,l},
\eeq
which is not identically zero in view of $\eqref{hypohn}$ and the fact that $H$ is generically of full rank.

On the other hand, for any $\b\in\N^d$, $\mdopz{t}{\b}P_{j+1}\left (\c^1,z^1,\c^2,\ldots;\Hb\rd\Vb_{j+1}(\c^1,z^1,\c^2,\ldots;t)\right )$ is convergent thanks to Lemma \ref{l:convQFSjHSj+1}. So, from Lemma \ref{l:lemdecv} we obtain the desired result.

\epf

We can now complete the proof of Theorem \ref{t:cvHgenecase}.

\bpf[Proof of Theorem {\rm\ref{t:cvHgenecase}}]

We use the previous setting and notation. From Proposition \ref{p:convHrdSj}, we have that for any $j\geq 1$, $H\rd V_{j}(\c^1,z^1,\c^2,\ldots;0)$ is convergent. Moreover, we known (see \cite[\S10.4, \S10.5]{BERbook}) that there exists $(p_k)_k$ a sequence which converges to 0 and such that, for $j=2(d+1)$, the rank of $V_{j}(\c^1,z^1,\c^2,\ldots;0)$ at any $p_k$ is $n+d$ and $V_j(p_k,0)=0$. So $V_{2(d+1)}(\c^1,z^1,\c^2,\ldots;0)$ has a convergent right inverse in a neighborhood of a suitable $p_k$ close enough to 0 which proves the convergence of $H$.
\epf

\subsection{Proof of  Theorem \ref{t:thcourbehol}}

To prove Theorem \ref{t:thcourbehol}, we need the following lemma.

\bl\label{l:lemcourbehol}
Let $(M,p)$ and $(\Mp,p^\p)$ be two germs of real-analytic hypersurfaces in $\C^{n+1}$. If $M$ is \hn \ at $p$  and if $\Mp$ does not contain any holomorphic curve through $p^\p$ then any holomorphic formal mapping sending $(M,p)$ into $(\Mp,p^\p)$ is either constant or generically of full rank.
\el

\bpf
The proof of Proposition 7.1 in \cite{BER00} contains the following fact: Let $M$ and $\Mp$ be formal real  hypersurfaces through points $p$ and $p^\p$ respectively and $H:(\C^{n+1},p)\rightarrow(\C^{n+1},p^\p)$ be a formal holomorphic  mapping sending $M$ into $\Mp$, with $M$ \hn \ at $p$ and $\Mp$ not containing any nontrivial formal holomorphic curve; if furthermore $M$ and $\Mp$ are given in normal coordinates and if $H=(F,G)$, then  $G\equiv0$ implies $H\equiv0$.

Moreover, in \cite{LM06} the authors proved, without any assumption on $\Mp$, that every formal holomorphic map $H:(\C^{n+1},p) \rightarrow(\C^{n+1},p^\p)$ sending $M$ into $\Mp$ which is transversally nonflat (i.e. satisfies in normal coordinates $G \nequiv 0$) satisfies $\J H \nequiv0$.

So, to prove the lemma it suffices to apply the two previous results noticing that if $\Mp$ is a real-analytic hypersurface that does not contain any holomorphic curve through $p^\p$ then it does not contain a formal holomorphic one through the same point (see \cite{Mil78}).
\epf

\bpf[Proof of Theorem {\rm\ref{t:thcourbehol}}]
The theorem is obtained by combining Lemma \ref{l:lemcourbehol} and Theorem \ref{t:cvHgenecase}.
\epf

\section{Proof of the Artinian Theorem}\label{s:pfartinresult}

\bpf[Proof of Theorem {\rm\ref{t:artinhyp}}] This proof is inspired by the proof of Lemma 14.2 in \cite{BMR02}. We use the notation of Section \ref{s:cvR} and in particular normal coordinates. From Theorem \ref{t:cvR} we know that $\mc{R}$ is convergent, thus using the Taylor expansion of $\Qp$, for any $\a\in\N^n$, $\Qpa(\Hb(\c,\t))=r_\a(\c,\t)$ converges. For fixed $k\in\N$, by Artin's approximation theorem \cite{Artin68} there exists $H^k:(\C^n\times\C,0)\rightarrow (\C^n\times\C,0)$ a holomorphic convergent power series mapping which agrees up with $H$ up to order $k$ at 0 and such that $\Qpa(\Hkb(\c,\t))=r_\a(\c,\t)$, for any $\a\in\N^n$. Consequently, we have the following power series identity,
\beq\label{egalQ}
\Qp(\zp,\Hb(\c,\t))=\Qp(\zp,\Hkb(\c,\t)).
\eeq

To show that, for every $k\in\N$, $H^k$ sends a neighborhood of 0 in $M$ into a neighborhood of 0 in $\Mp$, we consider 
\beq
\rp:(\zp,\wp,\cp,\tp)\in U\mapsto \wp-\Qp(\zp,\cp,\tp)\in\C,
\eeq
where $U$ is a sufficently small neighborhood of 0 in $\C^{2n+2}$, and 
\beq
\rpt(\zp,\wp,\cp,\tp):=\rpb(\cp,\tp,\zp,\wp).
\eeq
The convergent power series $\rpt$ is of rank 1 in a neighbourhood of 0 and $\rpt(\zp,\wp,\cp,\tp)=0$ implies $\rp(\zp,\wp,\cp,\tp)=0$ thanks to the mapping identity $\eqref{realcond}$. So, there exists $u\in\C[[\zp,\wp,\cp,\tp]]$ such that, 
\beq
\rp(\zp,\wp,\cp,\tp)=u(\zp,\wp,\cp,\tp).\rpt(\zp,\wp,\cp,\tp).
\eeq
We have, thanks to $\eqref{egalQ}$,
\beq\label{egalrp}
\rp\left (\zp,\wp,\Hb(\c,\t)\right )=\rp\left (\zp,\wp,\Hkb(\c,\t)\right ),
\eeq
and $H^k$ sends $M$ into $\Mp$ if and only if $\Qp\left (F^k(z,Q(z,\c,\t)\right ),\Hkb(\c,\t))=G^k\left (z,Q(z,\c,\t)\right )$ i.e. $\rp\left (H^k(z,Q(z,\c,\t)),\Hkb(\c,\t)\right )\equiv0$. But, 
\begin{displaymath}
\begin{array}{lll}
\rp\left (H^k(z,Q(z,\c,\t)\right ),\Hkb(\c,\t))&=&\rp\left (H^k(z,Q(z,\c,\t)),\Hb(\c,\t)\right ) \textrm{ from $\eqref{egalrp}$ }\\
&=&u\left (H^k(z,Q(z,\c,\t)),\Hb(\c,\t)\right ).\rpt\left (H^k(z,Q(z,\c,\t)),\Hb(\c,\t)\right ) \\
&=&u\left (H^k(z,Q(z,\c,\t)),\Hb(\c,\t)\right ).\rpb\left (\Hb(\c,\t),H^k(z,Q(z,\c,\t))\right ) 
\end{array}
\end{displaymath}

But, we have the identity $\rpb\left (\Hb(\c,\t),H^k(z,Q(z,\c,\t))\right )\equiv0$; Indeed, we recall that, since $H$ sends $M$ into $\Mp$, 
$$\Qp\left (F(z,Q(z,\c,\t)),\Hb(\c,\t)\right )=G(z,Q(z,\c,\t)).$$
So, applying $\Qpb\left(\Fb\left(\c,\t\right),F\left(z,Q\left(z,\c,\t\right)\right),.\right)$ to each side of the previous identity we obtain
\begin{displaymath}
\begin{array}{ll}
\Qpb\left(\Fb\left(\c,\t\right),F\left(z,Q\left(z,\c,\t\right)\right),\Qp\left(F\left(z,Q\left(z,\c,\t\right)\right),\Hb\left(\c,\t\right)\right)\right) \\
=\Qpb\left(\Fb\left(\c,\t\right),F\left(z,Q\left(z,\c,\t\right)\right),G\left(z,Q\left(z,\c,\t\right)\right)\right).
\end{array}
\end{displaymath}
Thus, using the mapping identity $\eqref{realcond}$, we have 
\beq
\Gb(\c,\t)=\Qpb\left (\Fb(\c,\t),H(z,Q(z,\c,\t))\right )
\eeq
i.e. $\rpb\left (\Hb(\c,\t),H(z,Q(z,\c,\t))\right )\equiv0$. Therefore $\rpb\left (\Hb(\c,\t),H^k(z,Q(z,\c,\t))\right )\equiv0$ and the proof of Theorem \ref{t:artinhyp} is complete.
\epf


\begin{thebibliography}{MMZ02a}

\bibitem[1]{Artin68}
  {\sc M. Artin} ---
  On the solutions of analytic equations.
  {\em Invent. Math.} {\bf 5} (1968), 277--291.

\bibitem[2]{BER97}
 {\sc M.S. Baouendi; P. Ebenfelt; L.P. Rothschild} ---
 Parametrization of local biholomorphisms of real analytic hypersurfaces.
 {\em Asian J. Math.} {\bf 1} (1997), 1--16.

\bibitem[3]{BERbook}
  {\sc M.S.~Baouendi; P.~Ebenfelt; L.P.~Rothschild} ---
  {\em Real Submanifolds in Complex Space and
  Their Mappings}. Princeton Math. Series {\bf 47}, 
  Princeton Univ. Press, 1999.

\bibitem[4]{BER00}
 {\sc M.S. Baouendi; P. Ebenfelt; L.P. Rothschild} ---
 Convergence and finite determination of formal CR mappings.
 {\em J. Amer. Math. Soc.} {\bf 13} (2000), 697--723.


\bibitem[5]{BMR02}
 {\sc M.S. Baouendi, N. Mir, L.P. Rothschild} ---
 Reflection ideals and mappings between generic submanifolds in complex space.
 {\em J. Geom. Anal.} {\bf 12} (4) (2002), 543--580.

\bibitem[6]{BM88}
 {\sc E. Bierstone, P.D. Milman} ---
 Semianalytic and subanalytic sets.
{\em Publ. Inst. Hautes \'{E}tudes Sci.  Math.} {\bf 67} (1988), 5--42.

\bibitem[7]{CM74}
 {\sc S.S. Chern, J.K. Moser} ---
 Real hypersurfaces in complex manifolds.
{\em Acta Math.} {\bf 133} (1974), 219--271.

\bibitem[8]{Gab73}
{\sc A.M. Gabri{\`e}lov} --- 
Formal relations among analytic functions.
{\em Izv. Akad. Nauk SSSR Ser. Mat.} {\bf 37} (1973), 1056--1090.

\bibitem[9]{Hua96}
  {\sc X. Huang} ---
  Schwarz reflection principle in complex spaces of dimension two. 
  {\em Comm. Partial Differential Equations} {\bf 21} (1996), 1781--1828.

\bibitem[10]{Juh08}
  {\sc R. Juhlin} ---
  Determination of formal CR mappings by a finite jet.
  To appear in {\em Adv. Math.} 2009. 

\bibitem[11]{LM06}
  {\sc B. Lamel, N. Mir} ---
  Remarks on the rank properties of formal CR maps.
  {\em Sci. China Ser. A} {\bf 49} (11) (2006), 1477--1490.


\bibitem[12]{MMZ03}
  {\sc F. Meylan, N. Mir, D. Zaitsev} ---
  On some rigidity properties of mappings between CR-submanifolds in complex space.
  {\em Journ\'ees \'Equations aux D\'eriv\'ees Partielles}, Exp.No. XII, 20 pp., Univ. Nantes, Nantes, 2003.

\bibitem[13]{Mil78}
  {\sc P. Milman} ---
  Complex analytic and formal solutions of real analytic equations in $\mathbb C^{n}$.
  {\em Math. Ann.}  {\bf 233} (1978), 1--7.

\bibitem[14]{Mir00} 
  {\sc N. Mir} ---
  Formal biholomorphic maps of real analytic  hypersurfaces.
  {\em Math. Res. Lett.} {\bf 7} (2000), 343--359. 

\bibitem[15]{Mir02}
  {\sc N. Mir} ---
  On the convergence of formal mappings
  {\em Comm. Anal. Geom.} {\bf 10} no. 1 (2002), 23--59. 

\bibitem[16]{Sta96} 
  {\sc N. Stanton} ---
  Infinitesimal CR automorphisms of real hypersurfaces.
  {\em Amer. J. Math.} {\bf 118} (1996), 209--233.   


 \end{thebibliography}
\end{document}